\documentclass[10pt]{article}
\usepackage{amsfonts,amssymb,latexsym, amsmath, amsthm}

\usepackage{amsfonts,amssymb,latexsym, amsmath, amsthm}

\usepackage[usenames]{color}%%  NEW

\theoremstyle{remark}

\def\ZZ{{\mathbb Z}}
\def\FF{{\mathbb F}}

\font\teneufm=eufm10
\font\seveneufm=eufm7
\font\fiveeufm=eufm5
\newfam\eufmfam
\textfont\eufmfam=\teneufm
\scriptfont\eufmfam=\seveneufm
\scriptscriptfont\eufmfam=\fiveeufm

\def\FF{{\mathbb F}}
\def\leq{\leqslant}
\def\geq{\geqslant}
\def\phi{\varphi}
\def\epsilon{\varepsilon}

\def\Aut{{\rm Aut}\,}
\def\Inn{{\rm Inn}\,}

\begin{document}

\def\normO{{\triangleleft_O\,}}

\def\C{{\rm C}}
\def\Z{{\rm Z}}
\def\cT{{\mathcal T}}
\def\cU{{\mathcal U}}

\def\wr{\,{\rm wr}\,}

\let\oldcomma=,
\catcode`\,=\active
\def,{\ifmmode \oldcomma \else {\/\rm\oldcomma}\fi}
\let\oldcolon=:
\catcode`\:=\active
\def:{\ifmmode \oldcolon \else {\/\rm\oldcolon}\fi}
\let\oldscolon=;
\catcode`\;=\active
\def;{\ifmmode \oldscolon \else {\/\rm\oldscolon}\fi}

\def\Cr{{\rm Cr}}
\centerline{PROFINITE GROUPS WITH FEW CONJUGACY CLASSES OF $p$-ELEMENTS}
\bigskip

\centerline{JOHN S.\ WILSON}

%\date{11 October 2021}
%\maketitle
\bigskip

{\small \noindent ABSTRACT. \
%\begin{abstract}
It is proved that a profinite group $G$ has fewer than $2^{\aleph_0}$ conjugacy classes of $p$-elements for an odd prime $p$ if and only if its $p$-Sylow subgroups are finite.  (Here, by a $p$-element one understands an element that either has $p$-power order or topologically generates a group isomorphic to $\ZZ_p$.) A weaker result is proved for $p=2$.}  %\end{abstract}
\bigskip

\centerline{1.\ INTRODUCTION}
\bigskip

Let $G$ be a profinite group and $p$ a prime.  We recall that a (generalized) $p$-element of $G$ is an element that (topologically) generates a finite $p$-group or a copy of $\ZZ_p$, and that the $p$-Sylow subgroups of $G$ are the maximal subgroups consisting of $p$-elements and are conjugate in $G$ (see \cite[(2.2.2)]{book}).   Consequently, if the $p$-Sylow subgroups of $G$ are finite, then $G$ has only finitely many conjugacy classes of $p$-elements.  

In \cite{jzn}, Jaikin-Zapirain and Nikolov proved that if $G$ is a profinite group with countably many conjugacy classes then $G$ is finite.  A small extension of their argument shows that if $G$ has fewer than $2^{\aleph_0}$ conjugacy classes then $G$ is finite.  Here we shall prove the following result.  

\begin{description}\item[Theorem A.] \begin{em}  Let $G$ be a profinite group and $p$ an odd prime.   If $G$ has fewer than $2^{\aleph_0}$ conjugacy classes of $p$-elements then its $p$-Sylow subgroups are finite.
\end{em}\end{description}

It follows at once that such a group $G$ has an open normal subgroup with trivial $p$-Sylow subgroup.

The corresponding result probably holds also for $p=2$, but here we only prove this for groups satisfying an extra condition.  We say that a finite simple group $S$ is $p$-{\em compliant} for a prime $p$ if either $p\nmid |S|$ or
for each $p$-element $h\in\Aut S\setminus \Inn S$ the coset $(\Inn S)h$ contains non-conjugate $p$-elements.   We shall prove (in Lemma 3.2 below) that all simple groups are $p$-compliant for $p$ odd.
It seems possible that they are $p$-compliant also for $p=2$, but checking this might be laborious. 
 
\begin{description}\item[Theorem B.] \begin{em} Let $G$ be a profinite group that has an open normal subgroup $N_1$ such that all non-abelian composition factors of finite continuous images of $N_1$ are $2$-compliant. If
$G$ has fewer than $2^{\aleph_0}$ conjugacy classes of $2$-elements then the $2$-Sylow subgroups of $G$ are finite. 
\end{em}\end{description}

{\hrule
\smallskip


\newpage

Let $G$ be a profinite group having non-trivial $p$-elements for infinitely many primes $p$.  
From a theorem of Herfort \cite{Wolfgang},  $G$ has a procyclic subgroup $A$ with the same property, and $A$ is isomorphic to a Cartesian product $A^*=\Cr_{p\in I}A_p$ of non-trivial procyclic $p$-groups $A_p$ for an infinite set $I$ of primes (see \cite[Chapter 2]{book}).  Elements of $A$ whose images in $A^*$ have distinct infinite supports have infinite order and generate non-isomorphic profinite subgroups, and certainly cannot be conjugate in $G$.

Thus a profinite group $G$ with fewer than $2^{\aleph_0}$ isomorphism classes of infinite procyclic subgroups has non-trivial $p$-Sylow subgroups for only finitely many primes; if 
$G$ also has fewer than $2^{\aleph_0}$ conjugacy classes of $p$-elements for each prime $p$ it follows from Theorems A and B, and the fact that each finite non-abelian simple group has order divisible by $3$ or $5$, that $G$ is finite. 
Therefore we recover the result of Jaikin-Zapirain and Nikolov \cite{jzn} mentioned above.

\bigskip

\centerline{2. \ PRELIMINARIES ON PROFINITE GROUPS}\medskip

For unexplained notation and general information about profinite groups and their Sylow theory, we refer the reader to \cite{book}.  We write $N\normO G$ to indicate that $N$ is an open normal subgroup of a profinite group $G$, and for
$x\in G$ we write $x^G$ for the conjugacy class of $x$ and $Nx^G$ for the product of the sets $N$, $x^G$.  Each conjugacy class $x^G$ is closed (being the image of $G$ under the continuous map $g\mapsto x^g$). 
All subgroups arising are understood to be closed.  A subset of a group is said to have {\em finite exponent} if its elements have bounded orders.   Here are some elementary facts that we shall use without further comment: each infinite profinite group has cardinality at least $2^{\aleph_0}$; each $p$-element in a quotient group lifts to a $p$-element; 
the property of having fewer than $2^{\aleph_0}$ conjugacy classes of $p$-elements is inherited by open subgroups and continuous images.  Our first lemma is well known.

\begin{description} \begin{em}\item[Lemma 2.1]  
Two elements $x,y$ of a profinite group $G$ are conjugate in $G$  if and only if for each $M \normO G$ the elements $Mx,My$ are conjugate in $G/M$. \end{em}
\end{description}

\begin{proof} \  If $Mx,My$ are conjugate for all $M$ then $y\in\bigcap_M Mx^G$, and this intersection equals $x^G$ as $x^G$ is closed (see \cite[Proposition 0.3.3]{book}).
\end{proof} 

We call an element of a profinite group a $p^\infty$-element if it generates (as a profinite group) a copy of $\ZZ_p$.

\begin{description} \begin{em}\item[Proposition 2.2.]  Let $G$ be a profinite group, $p$ a prime and $P$ a $p$-Sylow subgroup of $G$.  
Suppose that, for each coset $(N\cap P)t$ with $N\normO G$, $t\in P$, 

{\rm(a)} $(N\cap P)t$ contains elements not conjugate to $t$ in $G$, and

{\rm(b)}   $(N\cap P)t$ does not have finite exponent.

\noindent Then $G$ has at least $2^{\aleph_0}$ conjugacy classes of $p^\infty$-elements.
\end{em}\end{description}

\begin{proof}  We construct a descending chain $(N_k)_{k\geq1}$ of open normal subgroups and a family $(R_k)_{k\geq1}$ of finite subsets of
$P$ such that for each $x\in R_k$ there are elements $x^{(1)},x^{(2)}\in (N_{k-1}\cap P)x$ for which 

\ (i)  $N_k x^{(1)}$ and $N_k x^{(2)}$ are not conjugate in $G/N_k$ and

(ii) $N_kx^{(1)}$ and $N_kx^{(2)}$ have order at least $p^k$ in $G/N_k$.

\medskip

By (a) there is an element $w\in P\setminus\{1\}$.  Choose $N_1 \leq G$ with $w\not\in N_1$ and set $R_1=\{w\}$.  Clearly the required conditions hold.

Suppose that $k\geq2$ and that $N_{k-1}$, $R_{k-1}$ have been constructed, and let $x\in R_{k-1}$.  By (b) we can find $x^{(1)}\in (N_{k-1}\cap P)x$ and a subgroup $L_x\normO G$ with $L_x\leq N_{k-1}$ such that $L_xx^{(1)}$ has order at least $p^k$ in 
$G/L_x$.  By (a), Lemma 2.1 and again (b), we can find $M_x\normO G$ with $M_x\leq L_x$  and $x^{(2)} \in (N_{k-1}\cap P)x$ such that $M_xx^{(1)}$ and $M_xx^{(2)}$ are non-conjugate elements of $G/M_x$ and with $M_xx^{(2)}$ of order at least $p^k$ in 
$G/M_x$.  We define $N_k$ to be the intersection of the subgroups $M_x$
and set $R_k=\{x^{(1)},x^{(2)}\mid x\in R_{k-1}\}$.  

Now consider the family $F=\{ N_kx^G\mid k\geq1, x\in R_k\}$, partially ordered with respect to inclusion.  By construction, each element of $F$ contains two maximal elements; it follows that $F$ has
$2^{\aleph_0}$ maximal chains $(N_kx_k^G)$.  By compactness, the intersection $\bigcap N_kx_k^G$ of each such chain is non-empty; it is also evidently a union of conjugacy classes.  
If $(N_kx_k^G)$,  $(N_ky_k^G)$ are distinct chains then their intersections  $\bigcap (N_kx_k^G)$,  $\bigcap(N_ky_k^G)$ are disjoint: if $k$ is minimal with
 $N_kx_k^G\neq N_ky_k^G$, then $k\geq2$ and $x_k, ,y_k$ are distinct elements of $R_k$; since $N_{k-1}x_k^G=N_{k-1}y_k^G$ their relationship is that of the elements $x^{(1)},x^{(2)}$  appearing in the construction of $R_k$ and so $N_kx_k^G$ and $N_ky_k^G$ are disjoint. 
 
 Finally, if $x\in\bigcap (N_kx_k^G)$ then for each $r\geq 1$ we have $N_{r+1}x^{p^r}= N_{r+1}x_{r+1}^{p^r}\neq N_{r+1}$ and
 so $x^{p^r}\notin N_{r+1}$.  Therefore the elements of the intersections of chains $(N_kx_k^G)$ are $p^\infty$-elements and the result follows.
\end{proof}

\begin{description}\begin{em} \item[Lemma 2.3.] Let  $G$ be a profinite group, $p$ a prime and $P$ a $p$-Sylow subgroup of $G$.  
If each coset $(N\cap P)t$ with $N\normO G$ and $t\in P$
contains $p$-elements not conjugate to $t$ in $G$,
then $G$ has at least $2^{\aleph_0}$ conjugacy classes of $p$-elements. \end{em}
\end{description} 

\begin{proof} The statement is the same as that in Proposition 2.2, but without the hypothesis (b) on element orders and the conclusion that the $p$-elements obtained have infinite order.  The proof is also the same: we have $P\neq1$,  we construct sequences $(N_k)$, $(R_k)$ subject just to the condition (i) (with $x^{(1)}=x$ and $L_x=N_{k-1}$ for each $x$) and the proof is complete after the penultimate paragraph.
\end{proof}

 \begin{description} \item[Lemma 2.4.]  \begin{em} Suppose that $G$ is a pro-$p$ group with fewer than $2^{\aleph_0}$ conjugacy classes.  Then $G$ is finite.
 \end{em}\end{description}
 
 \begin{proof}  
 From Lemma 2.3 (with $G=P$) it is enough to consider the case when there exist $N\normO G$ and $t\in G$ such that $Nt\subseteq t^G$.  
 
Let $L\normO G$ with $L\leq N$ and write 
$\bar G$, $\bar N$, $\bar t$ for the images of $G, N$ and $ t$ in $G/L$.  Since $\bar N\bar t\subseteq \bar t^{\bar G}$ we have
$$|\bar N|=|\bar N\bar t|\leq |\bar G|/|\C_{\bar G}(\bar t)|,$$
and so $$|\C_{G/L}(Lt)|\leq |\bar G\colon \bar N|= |G\colon N|.$$
It follows from a result of Shalev \cite[Theorem A$'$]{shalev} that the derived length of $G/L$ is bounded in terms of $|G\colon N|$ alone.  Consequently $G$ is soluble.  If $G$ were infinite, it would have an open normal subgroup with an infinite abelian quotient group $A$, and $A$ too would have fewer than $2^{\aleph_0}$ conjugacy classes.  This is impossible and the result follows.  
\end{proof}

\begin{description}\item[Lemma 2.5.] \begin{em}  Let $G$ be a profinite group, $P$ a $p$-Sylow subgroup,  $N\normO G$, $t\in G$  and suppose that $(N\cap P)t$ has finite exponent.  Then  
$G$ has a finite series of closed characteristic subgroups in which each factor is one of the following: 
{\rm(i)} a pro-$p$-group; {\rm(ii)} a pro-$p'$-group, or {\rm(iii)} a Cartesian product of isomorphic finite simple groups of order divisible by $p$.  
\end{em}\end{description}

\begin{proof}  The result is very similar to \cite[Theorems 2 and 3]{oldcompact} and it is deduced from \cite[Theorems $2^\ast$ and $3^\ast$]{oldcompact} by exactly the same argument, described in the last paragraph of \cite[Section 2]{oldcompact}.
\end{proof}

We can now prove our theorems subject to an extra hypothesis (but without any hypothesis of $p$-compliance).

\begin{description}\item[Proposition 2.6.] \begin{em}  Let $G$ be a profinite group with fewer than $2^{\aleph_0}$ conjugacy classes of $p$-elements.  Let $P$ be a $p$-Sylow subgroup of $G$ and suppose that $(N\cap P)t$ has finite exponent for some open normal subgroup $N$ and some $t\in P$.  Then $P$ is finite.  \end{em}\end{description}

\begin{proof}  \
From Lemma 2.5, the group 
$G$ has a finite series of closed characteristic subgroups in which each factor is 
(i) a pro-$p$-group, (ii) a pro-$p'$-group, or (iii) a Cartesian product of isomorphic finite simple groups of order divisible by $p$.  By induction on the length of such a series, we can assume that $G$ has a normal subgroup $K$ of type (i), (ii) or (iii) with the property that the Sylow subgroup $PK/K$ of $G/K$ is finite.  

If $K$ is a pro-$p'$-group then the conclusion of the proposition clearly holds.  Suppose that $K$ is of type (iii).  Each $G$-invariant normal subgroup $L$ of $K$ has a $G$-invariant complement $\C_K(L)$, and so since there are maximal $G$-invariant subgroups of finite index it follows that $K$ is a Cartesian product of finite minimal $G$-invariant subgroups. The family of such minimal subgroups cannot be infinite, since $p$-elements chosen with supports that are unions of the elements of distinct subfamilies lie in distinct conjugacy classes.  
 Therefore $K$ is finite and the conclusion follows.

Suppose that $K$ is of type (i).  Since $PK/K$ is finite, $G$ has an open normal subgroup $M$ with $K\leq M\leq N$ such that $M/K$ and $PK/K$ are disjoint; the subgroup $M$ inherits the hypothesis on conjugacy classes and $M/K$ is a pro-$p'$-group. 
By the Schur--Zassenhaus theorem for profinite groups (cf.\ \cite[Proposition 2.33]{book}),  we have $M=K\rtimes Q$ for a pro-$p'$ subgroup $Q$.  We write $\Phi$ for the Frattini subgroup of $K$ (see \cite[Section 2.5]{book}) and regard $K/\Phi$ as an $\FF_pQ$-module.  Since $Q$ is a pro-$p'$-group, each finite continuous module image of $K/\Phi$ is a direct sum of simple modules by Maschke's theorem, and so $K/\Phi$ is a Cartesian product of simple modules.  However this Cartesian product cannot have infinitely many factors since elements with different supports lie in distinct $M$-orbits.  Thus $K/\Phi$ is finite, and an open subgroup $Q_0$ of $Q$ acts trivially on it.  It follows (for example from \cite[Theorem 5.1.4]{gor}) that $Q_0$ acts trivially on each finite continuous image of $K$ and hence acts trivially on $K$.  

Write $C=\C_Q(P)$. This subgroup is normal in $M$  and since $Q_0\leq C$ the image of $Q$ in $M/C$ is finite.  Therefore $KC/C$ is an open pro-$p$ subgroup of $M/C$ and it too inherits the hypothesis on conjugacy classes and so it is finite by Lemma 2.4.  Therefore the $p$-Sylow subgroups of both $M$ and $G$ are finite, as required. 
\end{proof}

\bigskip

\centerline{3. MINIMAL NORMAL SUBGROUPS OF FINITE GROUPS}

\medskip

Now we need to study conjugacy classes in cosets of minimal normal subgroups of a finite group $G$.  Our first lemmas elucidate the case when the minimal normal subgroup is a non-abelian simple group.

\begin{description}\item[Lemma 3.1.] \begin{em}  Let $S$ be a finite group with trivial centre, embedded naturally in its automorphism group $A$, and let $S\leq G\leq A$.  Let $p$ be prime with $p\mid|S|$ and $t$ a $p$-element of $A$.  Suppose that all $p$-elements of $St$ are conjugate to $t$ in $G$.  Then, passing to a subgroup of $G$ containing $S$, we can assume that $G$ has the following additional properties:

\noindent{\rm(a)} $[G,t]\leq S;$  
 
\noindent{\rm(b)} $G/S$ is a $p$-group$;$ and
  
\noindent {\rm(c) (i)}  if $p>2$ then $G/S$ is not metacyclic$;$ {\rm(ii)} if $p=2$ then $G/S$ is neither abelian nor a direct product of an abelian group and a dihedral group.
\end{em}\end{description}

\begin{proof}  Clearly $t\notin S$.  All relevant conjugations are done by elements of the preimage of $\C_{G/S}(St)$ and so we can replace $G$ by this preimage to assume (a).  Further, we replace $G$ by a subgroup minimal with respect to satisfying the hypothesis; in particular $$G=H\C_G(t) \hbox{ with }S\langle t\rangle\leq H\hbox{ implies }H=G.\eqno(\ast)$$

Each element $u$ of $St$ is uniquely a product of a $p$-element and a $p'$-element that commute with each other (they are both powers of $u$).  The $p'$-element must be in $S$, and so the $p$-element is in $St$ and is a conjugate of $t$.
Write $\C^-_S(x)$ for the set of $p'$-elements in $\C_S(x)$ for $x\in St$.  Each element of $\C^-_S(t^g)t^g$ is in $St$.  Therefore
$$|S|=|St|=|G/\C_G(t)||\C^-_S(t)|.$$  Set $k=|\C_S(t)|/|\C_S^-(t)|$.  Then 
$$|S|={{|G|}\over{k|\C_G(t)\colon \C_S(t)|}}={{|G|}\over{k|S\C_G(t)\colon S|}}$$ and so
$|G|=k|S\C_G(t)|$.  By Frobenius' Theorem \cite{frob}, $|\C_S^-(t)|$ is a multiple of the $p'$-part of $|\C_S(t)|$, and so $k$ must be a power of $p$.

Let $P$ be a Sylow $p$-subgroup of $G$ containing $t$.  Since $P$ and $S\C_G(t)$ have coprime indices we have $G=(PS)\C_G(t)$ and so $G=PS$ from $(\ast)$.  Assertion (b) follows.

Next, let $p>2$ and assume that $G/S$ is metacyclic.  Then $G/S$ is powerful, that is, its $p$-th powers generate its Frattini subgroup  $\Phi(G/S)$, as noted at the start of \cite[Section 3]{lumba}, and so for all $j\geq0$ the group generated by the $p^j$-th powers is the set of $p^j$-th powers, from \cite[Proposition 1.7]{lumba}.
Therefore for some $j\geq0$ the element $St$ of $G/S$ is the $p^j$-th power of an element $Sx$ not in $\Phi(G/S)$.  We may choose $x$ to be a $p$-element.   Thus $t$ is conjugate to $x^{p^j}$ and so there is $y\in G$ with $t=y^{p^j}$ and 
$Sy\notin \Phi(G/S)$.  If $G/S\neq \langle Sy\rangle$ then $G=H\langle y\rangle$ for a subgroup $H\geq S$ with $H<G$.  But since $y\in \C_G(t)$ we have $G=H\C_G(t)$, in contradiction to $(\ast)$.

Suppose now that $G/S= \langle Sy\rangle$; thus $G=S\langle y\rangle$ and since $y$ centralizes $t$ we have $G=S\langle t\rangle$ from $(\ast)$.  We use a generalization of Frobenius' Theorem  given in Hall \cite[Theorem III]{Hall}; this result implies that if $L$ is a group, $K\triangleleft L$ and $u\in L$, then for each integer $n\geq1$ the number of solutions of the equation $X^n=1$ in the coset $Ku$ is a multiple of the greatest common divisor $(|K|,n)$ of $|K|$ and $n$.  We take $L=G$, $K=S$, $u=t$ and let $n$ be maximum of $|\langle t\rangle|$ and $|S\cap P|$.  
Since $S$ must conjugate transitively the set $R$ of $p$-elements in $St$, this set is the set of solutions in $St$ of $X^n=1$, and $|R|$ is equal to $|S\colon \C_S(t)|$ and is a multiple of $(|S|,n)=|S\cap P|$; thus $\C_S(t)$ is a $p'$-group.  This is a contradiction since $t$ centralizes the intersection of $S$ and the centre of $P$.
 Assertion (c) (i) follows.  

Finally, if $p=2$ and $G/S$ is as in (i), then each element of $G/S$ is clearly a power of an element not in $\Phi(G/S)$ and we can argue as in the two paragraphs above. 
\end{proof}

The $p$-compliance of most simple groups now follows easily from known facts about simple groups.

\begin{description}\item[Lemma 3.2.] \begin{em} Let $S$ be a non-abelian finite simple group, embedded naturally in $\Aut S$.  Let $p$ be an odd prime divisor of $|S|$ and $t$ a $p$-element of $\Aut S$. 
Then the coset $St$ contains non-conjugate $p$-elements. 
\end{em}\end{description}

\begin{proof} We may assume that $St\neq S$.  We use the classification of the finite simple groups.  If $S$ is an alternating group $A_n$ with $n\neq6$ or a sporadic simple group then $|{\rm Out}\, S| \leq 2$ (see \cite[p.\ vii]{atlas}).  The remaining simple groups are Chevalley groups and twisted Chevalley groups, and  information about their outer automorphism groups is given in \cite[Table 5, p.\ xvi]{atlas}; see also
\cite[Theorem 2.5.12]{lyons}.  Inspection of this information shows that for $p\geq3$ 
the Sylow $p$-subgroups of ${\rm Out}\,S$ are always metacyclic. 
\end{proof}

We mention that $(\Aut S)/S$ has $2$-subgroups satisfying the conclusion of (ii) in Lemma 3.1 for some groups $S={\rm PSL}_n(q)$. 

We shall need the following elementary facts.

\begin{description}\item[Lemma 3.3.] \begin{em} Let $G$ be finite and $N\triangleleft G$.  Let $P$ be a Sylow $p$-subgroup of $G$ and $x\in P$.
If all elements in $(N\cap P)x$ are conjugate then so are
all $p$-elements in $Nx$. \end{em}\end{description}

\begin{proof} Since all relevant conjugations are done by elements of $\C_{G/N}(Nx)$ we may replace $G$ by the preimage of this subgroup.  Let $hx$ be a $p$-element of $Nx$. By Sylow's theorem it is conjugate to an element $h_1x$ of $P$, and then $h_1\in N\cap P$, so $h_1x$ is conjugate to $x$ by our hypothesis.  The result follows. \end{proof}

\begin{description}\item[Lemma 3.4.]\begin{em}  Let $U_1,\dots,U_n$ be non-empty sets such that for all $i,j$  the sets $U_i,U_j$ are either equal or disjoint. For $r=1,\dots,n$ let $V_r$ be a non-empty union
of some of the sets $U_i$.  If $$V_1\times\cdots \times V_n\subseteq \bigcup_{\sigma\in S_n} U_{1\sigma}\times\cdots\times U_{n\sigma}$$
 $($where $S_n$ is the symmetric group of degree $n)$ then $U_i=V_r$ for some $i,r$.
 \end{em}\end{description}
 
 \begin{proof}  The result is clear for $n=1$.  Let $n>1$.  
 After renumbering we can assume that $U_n\subseteq V_k$
 if and only if $k> m$.   Then $m<n$ since elements of 
 $V_1\times\cdots\times V_n$ have a co-ordinate in $U_n$.  
  For $u\in U_n$, we can have $(v_1,\dots,v_m,u,\dots, u)\in  U_{1\sigma}\times\cdots\times U_{n\sigma}$ only if $\sigma$ fixes $\{m+1,\dots,n\}$, and so
  $$V_1\times\cdots \times V_m\subseteq \bigcup_{\rho\in S_m} U_{1\rho}\times\cdots\times U_{m\rho}.$$
  The result follows by induction if $m\geq1$.  If $m=0$ then $U_n\subseteq V_r$ for all $r$.  
  
  Arguing similarly for each $U_i$ we obtain the result unless
  $\bigcup U_i\subseteq V_r$ for each $r$.   But then for $u\in\bigcup U_i$ the element $(u,\dots,u)$ is in some 
  $U_{1\sigma}\times\cdots\times U_{n\sigma}$; therefore $u$ is in each $U_i$ and all sets $U_i,V_r$ are equal. 
\end{proof}

Now we consider perfect minimal normal subgroups $M$ of finite groups $G$ in general.  Such a subgroup
$M$ is a direct product of all conjugates of a simple subgroup $S$ and it is helpful to pass to a permutational wreath product containing $G$.  The following lemma is very well known (see e.g.\ \cite[Lemma 3.1]{prob}).  

\begin{description}\item[Lemma 3.5.]\begin{em}Let $G$ be a finite group and $S$ a subgroup whose conjugates in $G$ generate their direct product $M$.
Let $\Gamma$ be the group of permutations of  the conjugates of $S$
induced by $G$ and
$A$ the image of ${\rm N}_G(S)$ in $\Aut S$, and suppose that $\C_G(M)=1$. Then there is an embedding of $G$ in $A\wr \Gamma$ in which $M$ maps onto the base group of $(\Inn S)\wr\Gamma$.
\end{em}\end{description}

To simplify calculations in permutational wreath products, 
we regard their base groups as groups of functions.  Thus if $H$ is a group and 
$\Omega$ a finite set we write $H^{(\Omega)}$ for the group 
of functions from $\Omega$ to $H$ with pointwise multiplication; if $\Gamma$ is a group of permutations of $\Omega$ we let $\Gamma$ act on the right on $H^{(\Omega)}$ as follows: 
$$f^\gamma(a)=f(a \gamma^{-1})\quad \hbox{for} \; 
f\in H^{(\Omega)}, a\in\Omega \hbox{ and } \gamma\in \Gamma.$$
The permutational wreath product 
$W=H\wr \Gamma$ is the split extension of $H^{(\Omega)}$ by $\Gamma$ with this action, and $H^{(\Omega)}$ is the base
group of $W$.  

\begin{description}\item[Lemma 3.6.]\begin{em} Let $M$ be a perfect minimal normal subgroup of a finite group $G$ and let $\C_G(M)=1$.  Thus $M$ is a direct product of conjugates of a simple group $S$.  Let $p$ be a prime that divides $|S|$ and $x$ a $p$-element in $G$.  If all $p$-elements of $Mx$ are conjugate in $G$ then there is a $p$-element
$u$ of $\Aut S$ such that all $p$-elements of $Su$ are conjugate in $\Aut S$.  
\end{em}\end{description}

\begin{proof}  Since conjugacy in $G$ implies conjugacy in a larger group, we can work with the wreath product of Lemma 3.5.  Let $H=\Aut S$, and $\Gamma ={\rm Sym}\, \Omega$; the minimal normal subgroup $M$ is now $S^{(\Omega)}$.  Let $P$ be a Sylow subgroup of $H$ and $Q=S\cap P$; then $P^{(\Omega)}$ and $Q^{(\Omega)}$ are $p$-Sylow subgroups of $H^{(\Omega)}$ and $M$ respectively.  Replacing $x$ by a conjugate we can assume that it lies in a Sylow subgroup of the wreath product whose intersection with $S^{(\Omega)}$ is
$Q^{(\Omega)}$; then every element of $Q^{(\Omega)}x$ is a $p$-element.  

Let $x=f\sigma$ with $f\in H^{(\Omega)}$, $\sigma\in \Gamma$.  First assume that $\sigma=1$.  The conjugate of $f$ by an element $d\rho$ is
$(d^{-1}fd)^\rho$.  For $\beta\in\Omega$ let $U_\beta$ be the $H$-conjugacy class of $f(\beta)$ and $V_\beta$
the set of $p$-elements in $Sf(\beta)$.  Thus $\prod V_\beta$ is the set of $p$-elements in $Mf$ and is contained in $f^G$, and 
the sets $U_\beta$, $V_\beta$ satisfy the same conditions as $U_i,V_i$ in Lemma 3.4.  Therefore
$U_\alpha=V_\beta$ for some $\alpha,\beta$ and our conclusion follows.

Next suppose $\sigma\neq1$.
Let $\beta\in \Omega$ lie in a $\langle\sigma\rangle$-orbit of length $r>1$ and define $u=\prod_{0\to r-1} f(\beta\sigma^i)$; the notation is intended to indicate that the factors in the product are taken with increasing $i$.

If $d\rho\in G$ and $d\rho$ conjugates $x$ to an element of $Mx$, then $\sigma, \rho$ commute and the conjugate is
$$d^{-\rho}f^\rho d^{\sigma^{-1}\rho}\sigma=(d^{-1}fd^{\sigma^{-1}})^\rho\sigma.$$ Suppose that $k\sigma\in Mx$ and $(k\sigma)^\rho$ is equal to the above element: then  for all $\beta\in\Omega$ we have
$d^{-1}(\beta)f(\beta)d(\beta\sigma)=k(\beta).$
Multiplying such equations over the $\langle\sigma\rangle$-orbit of $\beta$ we obtain
$$d^{-1}(\beta)u\,d(\beta)=d^{-1}(\beta)\bigg(\prod_{0\to r-1} f(\beta\sigma^i)\bigg)d(\beta)=\prod_{0\to r-1}k(\beta\sigma^i).$$

Now choose $l\in Q$ and define $k\in H^{(\Omega)}$ by $(kf^{-1})(\beta)=l$ and $(kf^{-1})(\gamma)=1$ for all $\gamma\neq\beta$.  Then $kf^{-1}\in Q^{(\Omega)}$ and $kf^{-1}x$ is a $p$-element in $Mx$.  Moreover
$\prod_{0\to r-1} (k(\beta\sigma^i))=lu$.  Defining $d$ as in the previous paragraph for this choice of $k$, we have
$d^{-1}(\beta)u\,d(\beta)=lu$.  Thus, from Lemma 3.4, every $p$-element of $Su$ is a conjugate in $H$ of $u$, as required.  
\end{proof}

The main result of this section now follows easily.

\begin{description}\item[Lemma 3.7.]\begin{em} Let $G$ be finite, $P$ a $p$-Sylow subgroup for a prime $p$, and $M$ be a perfect minimal normal subgroup with $p\mid|M|$.  
Assume in addition that the composition factors of $M$ are $p$-compliant. 
If $t$ is an element in $P\setminus M$ then $(M\cap P)t$ contains an element $y$ that is not conjugate to $t$.
\end{em}\end{description}
 
\begin{proof}  By Lemma 3.4 it suffices to show that $Mt$ contains such an element $y$, and  
this follows directly from the previous lemma, applied to $G/\C_G(M)$.    
\end{proof}

\bigskip

\centerline{4.\ PROOF OF THE THEOREMS}

\medskip

We still need some information about finite groups with no non-abelian composition factors of order divisible by  $p$.

\begin{description}\item[Lemma 4.1.]\begin{em}  Let $G$ be a  finite group and $p$ a prime. Let $N$ be a $p$-soluble normal subgroup and $x$ a $p$-element of $G$. If all
$p$-elements in $Nx$ are conjugate to $x$ then $N$ is an extension of a $p'$-group by a $p$-group.
\end{em}\end{description}

\begin{proof}  Suppose that the conclusion is false.  First we replace $N$ by its smallest normal subgroup of $p$-power index.
Next, since the hypothesis passes to quotients modulo normal subgroups in $N$, we can
suppose that $N$ is an extension of a minimal normal $p$-subgroup $A$ of $G$ by a $p'$-group; then
by the Schur--Zassenhaus Theorem $N=A \rtimes Q$ for a $p'$-subgroup $Q$ and by the Frattini argument $G=N\,D=AD$ where $D={\rm N}_G(Q)$.

Now $A\cap D\triangleleft G$, and if $A\cap D=A$ then $G={\rm N}_G(Q)$ and $N=A\times Q$,
a contradiction.  Therefore $G=A\rtimes D$.  As all elements of $Ax$ are conjugate, we may assume that
$x\in D$.  

For each $a\in A$ there exist $b\in A$, $d\in D$ with $ax=b^{-1}d^{-1}xdb$, hence with
$a=(b^{-1}b^{x^{-d}})(x^dx^{-1})$.  Since $G=A\rtimes D$ this equation implies that $x,d$ commute and that $a=b^{-1}b^x=[b,x^{-1}]$.
However $A\langle x\rangle$ is a finite $p$-group and so $[A,x^{-1}]<A$.  This is another contradiction.
\end{proof}

To prove our two theorems, because of Lemma 3.2 it suffices now to prove the following result.

\begin{description}\item[Lemma 4.2.]\begin{em}  Let $G$ be profinite, $p$ a prime, and suppose that $G$ has an open normal subgroup $N_1$ such that all non-abelian composition factors of finite continuous images of $N_1$ are $p$-compliant.  Let
$P$ be a $p$-Sylow subgroup of $G$. If $P$ is infinite then $G$ has at least $2^{\aleph_0}$ conjugacy classes of $p$-elements.
\end{em}\end{description} 

\begin{proof}   By Lemma 2.3 it suffices to prove that that each coset $(N\cap P)t$ with $N\triangleleft_O G$ and $t\in P$ contains non-conjugate elements.  
Let $N\triangleleft_O G$, $t\in P$.   

Find, if possible,  a subgroup $L\triangleleft_O G$ with $L<N\cap N_1$ such that $G/L$ has a perfect minimal normal subgroup $M/L\leq (N\cap N_1)/L$ of order divisible by $p$.  From Lemma 3.7, we can find an element $y\in(M\cap P)t$ such that $Lt,Ly$ are not conjugate in $G/L$; thus $y\in(N\cap P)t$ and $t,y$ are not conjugate in $G$.

Now suppose that no such $L$ can be found, and that all $p$-elements in $Nt$ are conjugate; then $N\cap N_1$ is a pro-($p$-soluble) group and applying Lemma 4.1 to each finite quotient we find that $N\cap N_1$ is an extension of a pro-$p'$ group by a pro-$p$ group $P_1$.  If $P$ is infinite then so is $P_1$, and it has at least $2^{\aleph_0}$ conjugacy classes by Lemma 2.5.  Therefore its preimage $N\cap N_1$ and $G$ both have $2^{\aleph_0}$ conjugacy classes of $p$-elements. The result follows. 
\end{proof}

Mathematisches Institut, Universit\"at Leipzig, 04109 Leipzig, Deutschland

\medskip

\quad and

\medskip

Christ's College, Cambridge CB2 3BU, United Kingdom

\medskip

E-Mail addresses:  {\tt wilson@math.uni-leipzig,de} and {\tt jsw13@cam.ac.uk} 

\end{document}